\documentclass[11pt]{article}
\usepackage{psfig,amssymb}

\newcommand{\df}{\displaystyle\frac}
\newcommand{\bsq}{\vrule height .9ex width .8ex depth -.1ex}

\newcommand{\hsp}{\hspace*{\parindent}}

\newcommand{\bmu}{{\mbox{\boldmath $\mu$}}}
\newcommand{\bpi}{{\mbox{\boldmath $\pi$}}}
\newcommand{\RR}{{\Bbb R}}

\newcommand{\NN}{{\Bbb N}}

\newcommand{\sL}{{\cal L}}
\newcommand{\sG}{{\cal G}}

\newcommand{\sM}{{\cal M}}

\newcommand{\bp}{{\bf p}}

\input amssym.def
\input amssym.tex

\newcommand{\beql}[1]{\begin{equation}\label{#1}}
\newcommand{\eeq}{\end{equation}}

\catcode`\@=11
\renewcommand{\section}{
        \setcounter{equation}{0}
        \@startsection {section}{1}{\z@}{-3.5ex plus -1ex minus
        -.2ex}{2.3ex plus .2ex}{\large\bf}
        }
\catcode`@=12
\makeatletter
\def\eqalignno#1{\displ@y \ta {\bf s} kip\@centering
  \halign to\displaywidth{\hfil$\@lign\displaystyle{##}$\ta {\bf s} kip\z@skip
    & $\@lign\displaystyle{{}##}$\hfil\ta {\bf s} kip\@centering
    & \llap{$\@lign##$}\ta {\bf s} kip\z@skip\crcr
    #1\crcr}}
\makeatother

\makeatletter
\def\@sect#1#2#3#4#5#6[#7]#8{\ifnum #2>\c@secnumdepth
     \def\@svsec{}\else 
     \refstepcounter{#1}\edef\@svsec{\csname the#1\endcsname.\hskip .75em }\fi
     \@tempskipa #5\relax
      \ifdim \@tempskipa>\z@ 
        \begingroup #6\relax
          \@hangfrom{\hskip #3\relax\@svsec}{\interlinepenalty \@M #8\par}%
        \endgroup
       \csname #1mark\endcsname{#7}\addcontentsline
         {toc}{#1}{\ifnum #2>\c@secnumdepth \else
                      \protect\numberline{\csname the#1\endcsname}\fi
                    #7}\else
        \def\@svsechd{#6\hskip #3\@svsec #8\csname #1mark\endcsname
                      {#7}\addcontentsline
                           {toc}{#1}{\ifnum #2>\c@secnumdepth \else
                             \protect\numberline{\csname the#1\endcsname}\fi
                       #7}}\fi
     \@xsect{#5}}
\def\@begintheorem#1#2{\it \trivlist \item[\hskip \labelsep{\bf #1\ #2.}]}
\makeatother
\begin{document}
\begin{center}
{\Large {\bf The Kruskal Count }} \\
\vspace{\baselineskip}
{\em Jeffrey C. Lagarias} \\
\vspace*{.2\baselineskip}
AT\&T Labs - Research \\
Florham Park, NJ 07932 \\

\vspace*{1\baselineskip}
{\em Eric Rains} \\
\vspace*{.2\baselineskip}
AT\&T Labs - Research \\
Florham Park, NJ 07932 \\

\vspace*{1\baselineskip}
{\em Robert J. Vanderbei} \\
\vspace*{.2\baselineskip}
Princeton University \\
Princeton, NJ 08544 \\

\vspace*{1\baselineskip}
(October 12, 2001) \\
\vspace{2\baselineskip}
{\em Abstract}
\end{center}

The Kruskal Count is a card trick invented by Martin Kruskal
in which a magician ``guesses'' a card selected by a subject
according to a certain counting procedure.
With high probability the magician can correctly ``guess'' the
card. The success of the trick is based on a mathematical principle
related to coupling models for Markov chains.
This paper analyzes in detail two simplified variants of the
trick and estimates the probability of success. The results
are compared with simulation data for several
variants of the actual trick.\\



\setlength{\baselineskip}{1.5\baselineskip}

%
%

\section{Introduction}
\hsp
The {\em Kruskal Count}  is a card trick invented by Martin D. Kruskal
(who is most well known for his work on solitons) which is described in
Fulves and Gardner \cite{FulGar75} and
Gardner \cite{Gar78},\cite{Gar88}.
In this card trick a magician ``guesses'' one card in a deck of cards which is
determined by a subject using a special counting procedure that we
call {\em Kruskal's counting procedure}.
The magician can with high probability identify
the correct card.

The subject shuffles a deck of cards as many times as he likes.
He mentally chooses a (secret) number between one and ten.
Kruskal's counting procedure then goes as follows.
The subject turns the cards of the deck face up one at a time,
slowly, and places them in a pile.
As he turns up each card he decreases his secret number by one
and he continues to count this way till he reaches zero.
The card just turned up at the point when the count reaches zero
 is called the
{\em first key card}
and its value is called the
{\em first key number.}
Here the value of an Ace is one, face cards are assigned the value five,
and all other cards take their numerical value.
The subject now starts the count over, using the first key number to determine
where to stop the count at the second key card.
He continues in this fashion,
obtaining successive key cards until the deck
is exhausted.
The last key card encountered, which we call the {\em tapped card}, is the
card to be ``guessed'' by the magician.

The Kruskal counting procedure for selecting the tapped card depends
on the subject's secret number and the
ordering of cards in the deck.
The ordering is known to the magician because the cards are turned
face up, but the subject's secret number is unknown.
It appears impossible for the magician to know the subject's
secret number. The mathematical
basis of the trick is that for most
orderings of the deck most secret numbers produce the same
tapped card.
For any given deck two different secret numbers produce two
different sequences of key cards, but if 
the two sequences ever have a key card in common, then they coincide
from that point on, and arrive at the same tapped card.
The magician therefore selects his own
secret number and carries out the Kruskal counting procedure for it
while the subject does his own count.
The magician's ``guess'' is his own tapped card.
The Kruskal Count trick succeeds with high probability, but if it fails the
magician must fall back on his own wits to entertain
the audience.

The problem of determining the probability of success of this trick leads
to some interesting mathematical questions.
We are concerned with the
{\em ensemble success probability}
averaged over all possible orderings of the deck (with the uniform
distribution).
Our objective in this paper is to estimate ensemble success probabilities for
mathematical idealizations of such counting procedures.
Then we numerically compare  the ensemble success
probabilities  on a 52-card deck with that of the Kruskal Count trick itself.
The success probability of the trick depends in part on the 
magician's strategy 
for choosing his own secret number.
We show that the magician does best to always choose the first card in the deck
as his first key card, i.e. to use secret number $1$.

The general mathematical problem we consider applies the Kruskal
counting procedure to a deck of $N$ labelled cards with each card label a
 positive integer, in which
each card has its label  drawn independently from some fixed
 probability distribution on the positive integers $\NN^+$. 
We call such distributions
{\em i.i.d. deck distributions};
they are specified by the probabilities $\{ \pi_j : j \geq 1 \}$ of a
fixed card having value $j$.
We assume that the subject  chooses an
initial secret number  from an initial  probability
distribution on $ \NN^+  =  \{ 1,~2,~3,~ \cdot \cdot \cdot  \} $,
and that the magician independently does the same from a
possibly different  initial probability distribution, and that
thereafter each follows the Kruskal counting procedure.
It is convenient to view the cards of the deck as turned over at unit times,
so that the card in the $M$-th position is turned over at time $M$.
If the $M$-th card is a key card for both magician and subject and no
previous card is a key card for both, then we say that
$M$ is the {\em coupling time}
for the sequences.
Let $t$ be a random variable denoting the coupling time on the
resulting
probability space
with $t  = +  \infty$ if coupling does not occur.
We wish to estimate  the ``failure probability''  ${\rm Prob}[t >  N]$.

The set of permutations of a fixed deck (with uniform distribution) does not
have the i.i.d. property, and is not Markovian, but  it
can be reasonably well
approximated by such a distribution.
The advantage of the simplifying assumption of an i.i.d. deck distribution
 is that the random variable $t$ can be interpreted as a
stopping time for a coupling method for a Markov
chain, as is explained in \S 2.

The mathematical contents of the paper are determination of
${\rm Prob}(t >  n)$ for a geometric i.i.d. deck distribution,
which is carried out in \S 3, and estimation
of $ {\rm Prob}(t >  n)$ for a uniform i.i.d. deck distribution,
which is carried out in \S 4. The proofs of several results 
stated in \S 4  are given in an appendix.

In \S 5 we consider the actual Kruskal count trick, and compare
its success probability with the approximations given by the
models  above. Because the Kruskal count trick using
an actual deck of 52~cards involves a stochastic process that
is not Markovian,
we estimate the success probability by Monte Carlo simulation. 
We consider the effect  on this  success
probability of varying the
magician's strategy for choosing his key card, and of varying the
value assigned to face cards. The magician should choose his
key card value to be 1. Assuming this strategy for the magician,
the success probability of the original Kruskal Count trick is
just over 85\%. 
Both the i.i.d. geometric distribution and i.i.d.
uniform distribution models above
give good approximations; the geometric distribution is off by
less than 3\%, and the uniform approximation is within 1\%.

There has been some previous work on mathematical
models of the  Kruskal count.
In 1975 Mallows \cite{Mal75} determined the expected value of
the coupling time of i.i.d. sequences, and observed that especially
simple formulae occur for the geometric distribution.
Recently Haga and Robins \cite{HR97} analyzed a simplified Markov
chain  model for the Kruskal count, which is related to, but
not the same as,  the models considered here. We discuss their
model further at the end of \S4. 

%
%

\section{Coupling Methods for Markov Chains}
\hsp
The
{\em coupling time
random variable} $t$ is a special case of a stopping time random
variable $t^\ast$ associated with a coupling method for studying a
Markov chain. This motivates our terminology.

To explain this connection, consider a homogeneous Markov chain
$(X_n : n \geq 0 )$ on a countable discrete state space $S$.
Given two initial probability distributions $\bp$ and $\bp'$
on $S$ a
{\em coupling method}
constructs a bivariate process $(X_n^1 ,~ X_n^2 )$
consisting of two copies of process $X_n$ with $X_0^1$ having
distribution $\bp$, $X_0^2$ having distribution $\bp'$,
and the two copies evolve independently until some (random)
{\em stopping time} $t^\ast$ at which
$X_{t^\ast}^1  = X_{t^\ast}^2$ and then
requires them to be equal thereafter, evolving as a single
process $X_n$. The stopping time $t^\ast$ is not necessarily required to be the
{\em first} time $t$ at which $X_t^1  = X_t^2$ occurs, and the particular
rule for choosing  $t^\ast$ defines the coupling method.
Let $ \bmu_n$, $\bmu'_n$ denote the distribution at time $n$
of the process $X_n$ stating from the distribution
$\bp,~\bp'$ respectively, at time 0, and let the
{\em variation distance}
$|| \bp - \bp' ||$ between two distributions on $S$ be
$$
|| \bp - \bp' ||  : = \sum_{s \in  S} | p(s) - p' (s) |~.
$$
The
{\em basic coupling inequality}
is
\beql{eq201}
\df{1}{2} || \bmu_n - \bmu'_n ||  \leq  {\rm Prob} [t^\ast >  n ]  ~.
\eeq
Such inequalities can be used to
prove ergodicity  of a Markov chain and to bound the speed of 
convergence to the
equilibrium distribution, by bounding the right side of
the inequality.

The first coupling method was invented by Doeblin \cite{Doe38}, and
many other coupling methods have been proposed since, see
Griffeath \cite{Gri78} for a survey.
Applications to card shuffling and random walks on groups are
described in Aldous and Diaconis \cite{AldDia86} and Diaconis \cite{Dia88}.
The basic coupling inequality (\ref{eq201})
is also valid for non-ergodic Markov chains, e.g. null-recurrent
or transient Markov chains on the state space $\NN$, as was observed by
Pitman \cite{Pit76}. 
Coupling
methods are traditionally used  as an auxiliary device to get 
information on the rate of convergence
to equilibrium of an ergodic Markov chain. 
In this paper,  
we are interested in obtaining upper and lower bounds for the coupling
probability itself, since it represents the failure probability of the
Kruskal Count trick. We do not
use the basic coupling inequality, but instead in
\S4  use inequalities relating coupling probabilities for
various different Markov chains.

For an i.i.d. deck the Kruskal counting procedure can be viewed as
moving on a Markov chain $\sM_\bpi$ on the state space $\NN$
where a state $j$ represents a current value of the
Kruskal counting procedure, with state 0 representing being at a key card, 
and state $j$ represents  that the next key card be
reached after exactly $j$ more cards are turned over. Each transition of
the Markov chain will correspond to turning over one card in
the deck.
Let the random variable $X_n$ denote the state of the Markov chain
at time n; it indicates the current Kruskal count value at
location $n$ of the deck, except that $X_n  = 0$ indicates a key
card at location $n$.
The transition probability for this chain from state 
$j \geq 1$ is probability 1 to state $j-1$
and 0 to all other states,
and from state 0 to state $j$ is
probability $\pi_{j+1}$, where $\{ \pi_j : j \geq 1 \}$ is
the distribution $\bpi$ of card labels. (That is, $\pi_1$ is the
probability that the key card has value $1$, and the chain 
transitions from state 0 to state 0.)
The initial distribution of secret numbers are distributions
$\bp ,~ \bp'$ on the state space $\NN$.
We define the random variable $t  = t( \bp ,~ \bp' )$ to be
the stopping time associated
to the coupling method that combines the chains $X_n^1$ and
$X_n^2$ at the first time that
$X_n^1  = X_n^2  = 0$. (This is not necessarily the first time
that $X_n^1  = X_n^2$.)
The basic coupling inequality (\ref{eq201}) for $\sM_\pi$ and $t$ then gives
\beql{eq202}
\df{1}{2} || \bmu_n - \bmu'_n || \leq {\rm Prob} [t >  n] ~,
\eeq
where $\bmu_n$ and $\bmu'_n$ are the $n$-step state probabilities 
for the chain $\sM_\pi$ started with initial distributions
$\bp$ and $\bp'$. 
We note that the Markov chain $\sM_\pi$ is ergodic if 
$E[ \bpi ]  = \sum_{j=1}^\infty j \pi_j$
is finite, and is null-recurrent otherwise.
In the ergodic case the stationary distribution
$\tilde{\bpi} = ( \tilde{\pi}_0 ,~ \tilde{\pi}_1 ,~ 
\tilde{\pi}_2 ,~ \cdot \cdot \cdot )$ is given by
$$
 \tilde{\pi}_j  = (1 - \pi_1 - \pi_2 - \cdot \cdot \cdot  - \pi_j )
(1  +  E[ \bpi ])^{-1}
$$
for $j \geq 0$.  This chain is ergodic for the deck distributions
that we consider, and our object is to estimate the ``failure probability''
${\rm Prob} [t >  n]$.

In the remainder of the paper, rather than considering Markov chains
of the type $M_{\bf \pi}$, we study simplified  Markov chains
that jump from one key card to the next, but which retain enough
information for coupling methods to apply.  

%
%

\section{Geometric Distribution}
\hsp
We consider an idealized deck consisting of cards  whose
labels are independently and identically distributed random
variables drawn from $ \NN^+ = \{ 1,~2,~3,  \cdot \cdot \cdot  \} $
with the geometric distribution $\sG_p$ given by $\pi_k  = (1  - p) p^{k-1}$,
$0 <  p < 1$. The geometric distribution has mean
\beql{eq30}
E[p]~~ =~~  \sum_{k=1}^{\infty}~  k \pi_k~~ =~~ \frac {1}{1 - p}.
\eeq
Let $\sG_N (p)$ denote the deck distribution induced on a deck  of $N$ cards.

Assume that the magician and subject both pick a secret number 
drawn from the same
geometric distribution $\sG_p$.
Let ${\rm  Prob} [  t >  N ] $ denote the probability
(choosing a deck of cards at random as above) that the
magician and subject have no common key card in positions 1 through $N$.

For the geometric deck distribution there is a simple 
exact formula for all coupling
probabilities. \\

\noindent
{\bf Theorem 3.1.} {\em For the geometric deck distribution $\sG_N (p)$ 
with initial geometric value distributions}  $\sG_p$,
\beql{eq301}
{\rm Prob}  [t >  N]   = p^N (2 - p)^N~.
\eeq
{\bf Proof.}   We use the
memorylessness property
of the geometric distribution,
which is that for a $\sG_p$-distributed variable $X$
the conditional probability ${\rm Prob} [u = k~ |~ u \geq \ell ]$ satisfies
\beql{eq302}
{\rm Prob} [u = k~ |~ u \geq \ell ]
= {\rm Prob} [u = k - \ell ] ~.
\eeq
By direct computation
$$
{\rm Prob}  [t >  1]  = 1 - (1  - p)^2  =  p(2 - p) ~.
$$
Now for $N \geq 2$,
\begin{eqnarray}
\label{eq303}
\mbox{Prob} [t > N] & = & \mbox{Prob} [t > N | X_1^1 \geq 2 ~{\rm and}~  
X_2^1 \geq 2 ]~~ \mbox{Prob} [X_1^1 \geq 2 ~{\rm and}~
X_2^1 \geq 2] \nonumber \\
  & & + \mbox{Prob} [t > N |X_1^1 = 1 ~{\rm and}~ X_1^2 \geq 2 ] ~~ \mbox{Prob} [X_1^1 = 1 ~{\rm and} ~X_2^1 \geq 2 ]  \\
  & & + \mbox{Prob} [t > N | ~X_1^1 \geq 2 ~{\rm and} ~X_1^2 = 1 ]  ~~ \mbox{Prob}[X_1^1 \geq 2  ~{\rm and} ~ X_2^1 = 1 ] \nonumber \\  
  & & + \mbox{Prob} [t > N | ~X_1^1 =1 ~{\rm and} ~X_1^2 = 1 ]  ~~ \mbox{Prob}
 [X_1^1 =1  ~{\rm and} ~ X_2^1 = 1 ], 
\nonumber
\end{eqnarray}
in which the last condition  $X_1^1 = X_2^1 = 1$ has zero 
probability for $N \geq 2$. 
Now by (\ref{eq302})
$$
\mbox{Prob} [t > N |X_1^1 \geq 2~{\rm and} ~X_1^2 \geq 2 ] =~
\mbox{Prob} [t > N -1 ] ~.
$$
In the second case $X_2^2 - 1$ is geometrically distributed, hence
by (\ref{eq302}) again
$$
{\rm Prob} [t > N | X_1^1 = 1 ~ \mbox{and}~
X_1^2 \geq 2 ] = {\rm Prob} [t > N -1 ] ~.
$$
The same holds for the third case, so (\ref{eq303}) becomes
\begin{eqnarray*}
{\rm Prob}  [t >  N]  & = &  {\rm Prob} [t >  N  - 1] {\rm Prob} 
[ \max (X_1^1 , X_1^2 ) \geq 2 ] \\
 &  =&  p(2  - p) {\rm Prob}  [t >  N - 1] ~.
\end{eqnarray*}
The theorem follows.  $~~~\bsq$

For the geometric distribution the magician can improve his chances by
 always selecting the first card.
Let $t'$ denote the coupling time for this process
where the subject draws his secret number from $\sG_p$.
Then one finds by a similar calculation that
\beql{eq304}
{\rm Prob}  [t' >  N]  = p(p(2  - p))^{N-1}= p^N(2 - p)^{N-1} ~,
\eeq
which is smaller than (\ref{eq301}) by a factor $1/(2-p)$.

%
%

\section{Uniform Distribution}
\hsp
Consider a deck of $N$ cards having a uniform i.i.d.
distribution of card values drawn from $[1,~B]$.
We estimate ${\rm Prob} [t >  N]$ where $t$ is the coupling time
assuming that both the magician and the subject draw a 
secret value {\em uniformly}
from $[1,~B]$.

For our analysis we introduce two auxiliary finite state Markov chains.
The first of these is a chain $\sL_B$ that we call the
{\em leapfrog chain.}
View the subject and magician as performing the Kruskal counting procedure
on two independently drawn decks.
The subject will use a white pebble to mark the location of key cards and the
magician will use a black pebble, according to their decks, and simultaneously
each moves to their respective first key card.
After this is done, the person having his pebble furthest
behind in the deck moves it to his next key card.
In case of a tie, where both pebbles are in the same relative position 
in the deck,
a move consists of both persons simultaneously moving their pebbles to
their next key cards, respectively. (Since the players have
separate decks, the next key card values of the two players
need not be the same.)
The states of the chain $\sL_B$ represent the distance the
white pebble is currently ahead of or behind the black pebble in the card
numbering, so there are $2B  -  1$ states $i$ with
$ - (B  - 1)  \leq  i  \leq  B - 1$.
A transition occurs whenever a pebble is moved; a transition from state
0 corresponds to both pebbles moving (independently), while a  transition
from any other state corresponds to exactly one pebble being moved.
A transition often involves one pebble leapfrogging
over the other, hence the choice of name for $\sL_B$.
The transition probabilities $p_{ij}$ are determined by the
uniform distribution on card values.
For $i  \neq  0$ the transition from $i$ to $j$
is determined by the value $v$ of the key card by
\beql{eq401}
v  = {\rm sign} (i) (i  - j) ~,
\eeq
so that
$$
p_{ij}  =  
 \left\{
\begin{array}{ll}
\df{1}{B} & \mbox{if ~ $1 \leq ~{\rm sign} (i) (i -j) \leq B$}~, \\
~~~ \\
0 & \mbox{otherwise} ~,
\end{array}
\right. \eqno{\mbox{(4.2a)}}
$$
while for $i = 0$ the transition probabilities are
$$
\pi_j  : =  p_{0j}  = \df{B - |j|}{B^2} ~. \eqno{\mbox{(4.2b)}}
$$
This chain is ergodic, and it is easy to check that $\pi_j$
in (4.2b) gives the stationary distribution for $\sL_B$.
Table 4.1 gives the state transition matrix
$[ p_{ij} ]$ for $\sL_4$.
\begin{table}
\begin{center}
\begin{tabular}{|cc|c|c|c|c|c|c|c|} \hline
 & & \multicolumn{7}{c|}{exit state $j$} \\ \cline{3-9}
  & & $-3$ & $-2$ & $-1$ & 0 & 1 & 2 & 3 \\ \hline
 & 3 & 0 & 0 & 1/4 & 1/4 & 1/4 & 1/4 & 0 \\
 & 2 & 0 & 1/4 & 1/4 & 1/4 & 1/4 & 0 & 0 \\
 entering & 1 & 1/4 & 1/4 & 1/4 & 1/4 & 0 & 0 & 0 \\
 state & 0 & 1/16 & 2/16 & 3/16 & 1/4 & 3/16 & 2/16 & 1/16 \\
 $i$ & -1 & 0 & 0 & 0 & 1/4 & 1/4 & 1/4 & 1/4 \\
 & -2 & 0 & 0 & 1/4 & 1/4 & 1/4 & 1/4 & 0 \\
 & -3 & 0 & 1/4 & 1/4 & 1/4 & 1/4 & 0 & 0 \\ \hline
\end{tabular}
\end{center}
\begin{center}
Table 4.1:Leapfrog Chain $\sL_4$  
\end{center}
\end{table}

Now consider the case that the subject and magician perform the 
Kruskal counting procedure on the {\em  same} deck.
As long as their sequences of key cards remain disjoint, these 
key card values are independent
random variables, and their relative positions of current key cards
are described by transitions of the leapfrog chain.
This persists until they have a
key card in common, i.e. until the
state 0 is reached on the leapfrog chain.
Thus ${\rm Prob} [t >  N]$ corresponds to the probability of those
sequences of transitions in the leapfrog chain starting from 0
that avoid the 0 state until one pebble has moved to a position
beyond $N$.
We can keep track of sequences that never visit 0 by forming the
{\em reduced leapfrog chain}
$\bar{\sL}_B$ obtained by deleting the 0 state and assigning new 
transition probabilities
\setcounter{equation}{2}
\beql{eq403}
\bar{P}_{ij} := (1- p_{i0} )^{-1} p_{ij}~.
\eeq
For $\sL_B$ the probability of
going to 0 is a constant, hence
\beql{eq404}
\bar{P}_{ij} = \left( 1 - \df{1}{B} \right)^{-1} p_{ij} ~,
\eeq
so that all
values $\bar{P}_{ij}$ are either $\frac{1}{B  - 1}$ or 0.
Table 4.2 gives the state transition probabilities
$[ P_{ij}]$ for $ \bar{\sL}_4$.
\begin{table}
\begin{center}
\begin{tabular}{|cc|c|c|c|c|c|c|} \hline
 & & \multicolumn{6}{c|}{exit state $j$} \\ \cline{3-8}
  & & $-3$ & $-2$ & $-1$  & 1 & 2 & 3 \\ \hline
 & 3 & 0 & 0 & 1/3 & 1/3 & 1/3  & 0 \\
 & 2 & 0 & 1/3 & 1/3 & 1/3  & 0 & 0 \\
 entering & 1 & 1/3 & 1/3 & 1/3  & 0 & 0 & 0 \\
 state & -1 & 0 & 0 & 0 & 1/3 & 1/3 & 1/3 \\
 $i$ & -2 & 0 & 0 & 1/3 & 1/3 & 1/3 & 0 \\
 & -3 & 0 & 1/3 & 1/3 & 1/3 & 0 & 0 \\ \hline
\end{tabular}
\end{center}
\begin{center}
 Table 4.2: Reduced Leapfrog Chain $\bar{\sL}_4$
\end{center}
\end{table}

The
{\em initial state distribution}
on the reduced leapfrog chain $\bar{\sL}_B$ corresponds to that
after one transition of the leapfrog chain from the 0 state,
conditioned on not staying at 0, which is
\beql{eq405}
\bar{\pi}_j  : =  \left(  1 -  \df{1}{B} \right)^{-1}
\df{B - |j|}{B^2} ~, ~~~ \mbox{for 
$1  \leq  |j|  \leq  B  - 1$} ~.
\eeq
This chain is ergodic and has $\bar{\pi}_j$ as its stationary
distribution.

We next define a random variable
$t_{N,B}^{\ast \ast}$
which counts the total number of key
cards produced during the Kruskal count by the subject and
magician, up to and including the first key card that occupies
 a position exceeding $N$. We call $t_{N,B}^{\ast \ast}$
the {\em travel time beyond position N}.
To determine
the travel time, we require as additional  data the
position $i$ of the
{\em top key card,}
which we define to be that key card which is closest to the
top of the deck.
Given that the initial state of the chain is in state $j$
the conditional probability $r_{ij}$ that the top key
card is in position $i$ is
\beql{eq406}
r_{ij}  = \df{1}{B  - |j|} ~~~~~ 1  \leq  i  \leq  B  - |j| ~,
\eeq
and is 0 otherwise.
The position of the top key card together with the sequences
of successive states of $\bar{\sL}_B$ allow the reconstruction of
all moves during the Kruskal count, and the determination of
the travel time $t_{N,B}^{\ast \ast}$. \\
{\bf Lemma 4.1. }
{\em If $N \geq B$ then}
\beql{eq407}
{\rm Prob}  [t >  N]  = \sum_{j=1}^N \left(  1 - \df{1}{B} \right)^{j-1} 
{\rm Prob} [ t_{N,B}^{\ast \ast}  = j] ~.
\eeq

{\bf Proof. }
The event $[t >  N]$ corresponds to all sequences of state
transitions in $\bar{\sL}_B$  starting at state 0 that never
return to 0 before some pebble moves to a position $ \geq N  +  1$.
Such a sequence of transitions is matched (after the first move)
by corresponding state transitions in $\bar{\sL}_B$.
The probabilities between $\sL_B$ and $\bar{\sL}_B$
differ by a multiplicative factor $\left(  1 - \frac{1}{B} \right)$.
There is one less factor of $\left(  1 - \frac{1}{B} \right)$
than $t_{N,B}^{\ast \ast}$ counts because the initial state
of $\bar{\sL}_B$ counts as two key cards, but corresponds to
only one transition in $\sL_B$.  $~~~\bsq$

Lemma 4.1 is useful because the distribution of the travel time
$t_{N,B}^{\ast \ast}$ is strongly peaked and
relatively tractable to estimate.
Since no move of a pebble is larger than $B$, and since
both pebbles are within $B$ cards of the $N$-th card at the
stopping time $t_{N,B}^{\ast \ast}$, one has
\beql{eq408}
t_{N,B}^{\ast \ast} \geq \df{2N}{B} - 1 ~.
\eeq
Lemma 4.1 then yields
\beql{eq409}
{\rm Prob} [t >  N]  \leq  \left(  1 - \df{1}{B} \right)^{\frac{2N}{B}  - 2} ~.
\eeq
This shows the (well-known) fact that ${\rm Prob} [t >  N]$ decreases
exponentially as a function of $N$.

Using large-deviation theory we can obtain the asymptotic behavior of
${\rm Prob}  [t >  N]$ as $N \rightarrow \infty$. \\
{\bf Theorem 4.1. }
{\em For fixed $B$ there is a  positive constant $\alpha_B$ such that}
\beql{eq410}
{\rm Prob} [t >  N]  = \exp ({ -\alpha_B N (1  +  o(1))} )
\eeq
{\em as $N \rightarrow \infty$.}

We relegate the proof of this result to the appendix,  where
 we also give a variational formula for $\alpha_B$.
We easily obtain from (\ref{eq409}) the inequality 
\beql{eq410aa}
\alpha_B \geq  ( \frac{2}{B}) | \log ( 1 - \frac {1}{B}) |  = 
\frac{2}{B^2} + O(\frac {1}{B^3}).
\eeq  
It is intuitively clear that the  expected value of a key card is
$ \geq  \frac{B}{2}$ in all states, hence one expects that
${\rm Prob} \left[ t_{N,B}^{\ast \ast}  \leq  \frac{N}{B} \right]  \geq  \frac{1}{2}$,
which with Lemma 4.1 would imply that $\alpha_B \leq \frac{4}{B^2} 
+ O(\frac {1}{B^3})$.
Theorem 4.2 below shows that $\frac {B^2}{\alpha_B} \rightarrow 4$
as $B \rightarrow \infty$, see (\ref{eq419aa}).

We  next obtain upper and lower bounds for 
${\rm Prob}  [t_{N,B}^{\ast \ast} >  k]$
by approximating the reduced leapfrog chain $\bar{\sL}_B$ with two
simpler Markov chains
$\sL_B ^+$ and $\sL_B^-$, as follows.
These chains both describe the leapfrog motion of two colored pebbles
at most $B$ units apart, with the states representing the current
distance the white pebble is ahead.
\begin{description}
\item{\hspace{.5in}(i)}
 In $\sL_B^-$ the
pebble further behind jumps $v$ units with $v$ drawn uniformly
from the range $[1,~B  - 1]$.
\item{\hspace{.5in}(ii)}
In $\sL_B^+$ the pebble further behind jumps $v$ units
with $v$ drawn uniformly from the range $[2,~B]$.
\end{description}
The chain $\sL_B^-$ is exactly like the leapfrog chain $\sL_{B-1}$
except that in state 0 only the white pebble jumps.
The chain $\sL_B^-$ has $2B  - 1$ states labelled by
$|i|  \leq  B  - 1$, while the chain
$\sL_B^+$ has $2B  +  1$ states labelled by $|i|  \leq  B$.
Both chains have the property that the card values drawn are
independent of the current state.
For the chain $\sL_B^-$ we define a {\em travel time} 
$t_{N,B}^-$  beyond position N, which is obtained by starting 
the chain in state 0,
with both pebbles in position 0,
associating a movement of pebbles on a line with each
state transition, and counting the total number of state transitions
up to and including the first time that a pebble is moved beyond
position $N$.
For the chain $\sL_B^+$ we define a {\em travel time}
$t_{N,B}^+$  beyond position N similarly.\\
{\bf Lemma 4.2. }
{\em For all $N,~B$ and $k$, one has}
\beql{eq411}
{\rm Prob} [t_{N,B}^+ > k ]~ \geq~~{\rm Prob}
[t_{N,B}^{\ast \ast} > k ]~ \geq~ {\rm Prob}
[t_{N,B}^- > k ] ~.
\eeq

We give the proof of Lemma 4.2 in the appendix.
Lemmas 4.1 and 4.2 when combined yield the bounds
\beql{eq415}
P_{N,B}^+ \geq {\rm Prob} [t > N ] \geq P_{N,B}^-
\eeq
where
$$
P_{N,B}^\pm := \sum_{j=1}^N \left( 1 - \df{1}{B} \right)^{j-1}
{\rm Prob} [t_{N,B}^\pm =  j ] ~.
$$
The simple form of the chains $\sL_B^+$ and $\sL_B^-$
allows the asymptotic behavior of $P_{N,B}^+$ and
$P_{N,B}^-$ to be explicitly determined, as follows. \\
{\bf Theorem 4.2. }
{\em For fixed $B$ as $N \rightarrow \infty$ one has}
\beql{eq416}
P_{N,B}^\pm   = \exp (-   \alpha_B^\pm N (1  +  o(1)))
\eeq
{\em where $\frac{1}{2} \alpha_B^-$ is the unique root $\alpha$ of}
\beql{eq417a}
\sum_{i=1}^{B-1}  \exp (i \alpha )  = B~, 
\eeq
{\em and $\frac {1}{2} \alpha_B^+$ is the unique root of}
\beql{eq417b}
\sum _{i=1}^{B-1} \exp  ((i  +  1) \alpha )  = B~. 
\eeq
{\em As $B  \rightarrow \infty$ these quantities satisfy}
\beql{eq418a}
\alpha_B^+   =   \df{4}{B^2 } - \df{20/3}{B^3}  + 
O(B^{-4} ) 
\eeq
\beql{eq418b}
\alpha_B^-  =  \df{4}{B^2}  + \df{4/3}{B^3}  + 
O(B^{-4}). 
\eeq
\\

The proof of this result is given in the appendix.
Theorem 4.1 together with the inequalities (\ref{eq415}) shows that 
for large $B$
one has
\beql{eq419aa}
{\rm Prob}[t >  N]  = \exp \left(   - \left( 
\df{4}{B^2}  +  O \left( \df{1}{B^3 } \right) \right)
(1  +  o (1) ) N \right)
\eeq
as $N \rightarrow \infty$. 

We relate these results to the model of Haga and Robins \cite{HR97}.
The Markov chain studied by  Haga and Robins is obtained
from the leapfrog chain by identifying states $k$  and $-k$ for
all $k \geq 1$; thus it has exactly $B$ states.
The resulting chain factors out the action of the involution sending
$k$ to $-k$ under which the chain probabilities are invariant,
and this loses the ``leapfrog'' information which is
necessary for computing  exact coupling probabilities.
Haga and Robins estimate instead the
probability of avoiding absorption in the absorbing state $0$ in
the first M transitions of the resulting factor chain.
This probability asymptotically decays like $O( (\lambda_B)^M)$
as $M \to \infty$, where
$\lambda_B$ is 
the modulus of the second
largest eigenvalue 
 of the characteristic polynomial of their Markov chain. 
The characteristic polynomial of the transition matrix
of the Haga-Robins Markov
chain is
$p_B(x) := (x + 1/B)^B - (1 + 1/B)^B x^{B-1},$
and it  can be shown that the modulus of its second largest eigenvalue 
satisfies
\beql{eq411b}
\lambda_B = 1 - \frac{2}{B} + O(\frac{1}{B^2})
\eeq
as $B \to \infty$.
To relate  $\lambda_B$ to 
the asymptotic coupling probability  decay rate $exp(-\alpha_B)$
in Theorem 4.1,
we note that the expected size of a step in the Haga-Robins
chain is about $B/2$, so that after $M$ steps the location of
the chain should be around the position
$N \approx MB/2$.
One should therefore compare  $(\lambda_B)^{2/B} $ and $exp(-\alpha_B)$,
and one finds 
that both of these
quantities 
are asymptotic to $1 - \frac{4}{B^2} + O(\frac{1}{B^3})$
as $B \to \infty,$  using (\ref{eq411b}) and Theorem 4.2.

%
%
\newpage
\section{Numerical Results: The Kruskal Count}

We compare predictions obtained from the two models studied in this
paper with the 
performance of the actual Kruskal Count trick. 

For the actual
Kruskal count we consider a standard deck of 52 cards, and
we assume that the subject draws a key card using a uniform
distribution from the set of available key card values. We study
the effects of varying the magician's strategy  on the success
probability of the Kruskal
Count trick. The magician has the freedom to choose his key
card, and he has also has the extra freedom to  specify a rule
for assigning values to the ``face cards'' {J, Q, K}. We
study three possible rules variations:
\begin{description}
\item{\hspace{.5in}(a)}
Assign the values {11, 12, 13} to {J, Q, K}, respectively.
\item{\hspace{.5in}(b)}
Assign the value 10 to  each of {J, Q, K},
\item{\hspace{.5in}(c)}
Assign the value 5 to each of {J, Q, K}.
\end{description}
The first two of these rules variations are presented as ``straw men''
useful for comparison with the models of this paper.  To obtain numerical
values for the Kruskal count trick we used a Monte Carlo simulation with
$10^6$ trials for each data point.  For simulations of the i.i.d. uniform
deck distribution, an ``exact'' calculation was done using an enlarged
Markov chain which kept a running total of the value of the position N of
the leading pebble, and with an absorbing state whenever a pebble jumps
past the end of the deck.  Since the smallest step size is 1, this chain
reaches an absorbing state after a number of steps equal to the size of the
deck; consequently, we it suffices to compute the state of the chain after
that number of steps.  Simulations of the i.i.d. ``semiuniform''
distributions for rules variations $(b)$ and $(c)$ were done similarly to
the i.i.d. uniform case.

The rules variation $(a)$ corresponds to the uniform distribution
on $\{ 1,2,...,13 \}$. The average key card size is $7$.  We therefore
consider as an approximation the i.i.d. geometric deck distribution
with $p = \frac {6}{7}$, which has mean key card size $7$.  According
to Theorem 3.1 the failure probability for$N=52$ the magician drawing 
his first key card according to the geometric distribution is
\beql{eq501}
FP(\sG_a) = (\frac {6}{7})^{52} (\frac {8}{7})^{52} = 0.342254.
\eeq
If the magician chooses the first card to be his first key card,
by (3.5) his failure probability for $N=52$ is
\beql{eq502}
FP({\sG'}_a) = \frac{7}{8} FP(\sG_a) = 0.299472
\eeq

\begin{table}
\begin{center}
\begin{tabular}{|c|c|c|}
\hline
&Kruskal&uniform\\
\hline
1  & 0.315180 & 0.319486 \\
2  & 0.318564 & 0.322994 \\
3  & 0.321975 & 0.326492 \\
4  & 0.325298 & 0.329981 \\
5  & 0.328794 & 0.333461 \\
6  & 0.332235 & 0.336929 \\
7  & 0.336055 & 0.340385 \\
8  & 0.339264 & 0.343827 \\
9  & 0.342758 & 0.347251 \\
10 & 0.346780 & 0.350655 \\
11 & 0.349464 & 0.354034 \\
12 & 0.353367 & 0.357385 \\
13 & 0.357044 & 0.360703 \\
\hline
avg & 0.335906 & 0.340276 \\
\hline
\end{tabular}
\end{center}
\begin{center}
Table 5.1: Failure probabilities for rules variation $(a)$
\end{center}
\end{table}

Table~5.1 presents data for rules variation $(a)$ for the Kruskal
Count and the i.i.d. uniform deck distribution on $\{ 1,2,...,13 \}$.
The table gives failure probabilities in which the magician's 
strategy is to choose as first key card the j-th card, for 
$1 \leq j \leq 13$. plus a final row that gives the failure
probability when the magician draws a card uniformly in $\{ 1,2,...13 \}$.
The data in Table~5.1 show that the magician does best to choose
$j=1$ as his key card. The non-Markovian nature of the actual deck
causes the failure probabilities to differ from the i.i.d. uniform
deck distribution; the effect is a decrease of about 0.3\%.
We also see that the failure probability for the i.i.d. geometric 
distribution  is an overestimate of
the failure probability for the Kruskal Count when the magician 
picks a random card
as first key card,
and underestimates the failure probability when the magician
picks the first card as key card. 

We next consider  the rules variations $(b)$ and $(c)$. 
For rules variation
$(b)$ the expected key card size is $\frac {85}{13}$, so for comparison
we consider the
i.i.d. geometric deck distribution with $p = \frac{72}{85}$. If the
magician chooses his first key card according to the same geometric 
distribution, then the failure probability is
\beql{eq503}
FP(\sG_b) = (\frac {72}{85})^{52}(\frac {98}{85})^{52}
= 0.292064,
\eeq
while if the magician draws the first card as his key card, then
\beql{eq504}
FP({\sG'}_b )~ =~\frac {85}{98} FP(\sG_b) = 0.253320.
\eeq
For rules variation $(c)$ the expected key card size is $\frac{70}{13}$,
so for comparison we consider the i.i.d. geometric deck distribution
with $p = \frac {57}{70}$. If the magician chooses his first key card
with the same geometric distribution, then the failure probability is
\beql{eq505}
FP(\sG_c) = (\frac {57}{70})^{52} (\frac {83}{70})^{52} = 0.161197,
\eeq
while if the magician chooses the first card as his key card, the
failure probability is
\beql{eq506}
FP({\sG'}_c ) = \frac {70}{83} FP(\sG_c) = 0.135949.
\eeq

\begin{table}
\begin{center}
\begin{tabular}{|c|c|c||c|c||c|}
\hline
& Kruskal (b) &   semiuniform (b)&   Kruskal (c) &   semiuniform (c) & uniform\\
\hline
1  & 0.277869 & 0.284060 & 0.146238 & 0.152658 & 0.150944\\
2  & 0.280756 & 0.287235 & 0.148801 & 0.155266 & 0.153684\\
3  & 0.284330 & 0.290447 & 0.151204 & 0.157847 & 0.156407\\
4  & 0.287163 & 0.293623 & 0.153736 & 0.160399 & 0.159109\\
5  & 0.290317 & 0.296782 & 0.156075 & 0.162918 & 0.161789\\
6  & 0.293557 & 0.299920 & 0.159744 & 0.166357 & 0.164444\\
7  & 0.296910 & 0.303034 & 0.162474 & 0.168973 & 0.167070\\
8  & 0.300023 & 0.306118 & 0.164977 & 0.171553 & 0.169665\\
9  & 0.303194 & 0.309171 & 0.167735 & 0.174094 & 0.172225\\
10 & 0.306383 & 0.312185 & 0.170064 & 0.176591 & 0.174747\\
\hline
avg& 0.292050 & 0.298258 & 0.158105 & 0.164666 & 0.163008\\
\hline
\end{tabular}
\end{center}
\begin{center}
Table 5.2: Failure probabilities for rules variations (b) and (c)
\end{center}
\end{table}

Table~5.2 presents failure probability data for rules variations 
$(b)$ and $(c)$ for the
Kruskal Count and for the i.i.d. semiuniform  deck distributions 
which have the
card values $\{ 1,2,...,10 \}$ chosen with the same probabilities as 
rules variations $(b)$ and
$(c)$ impose on the actual deck. The non-Markovian nature 
of the actual deck results in the Kruskal count
failure probabilities differing from the  corresponding i.i.d. 
deck distributions; they are smaller by about 0.6\%. 
The failure probability for the i.i.d. geometric
distribution when the magician chooses the first card as first key card
gives an underestimate for the failure probabilities
of the Kruskal Count in rules variations $(b)$ and $(c)$. 
The numerical results show that
the magician should choose the first card as his key card. The
effect of the choice of the magician's key card on
the failure probability is small, at most 2.5\%.
 In comparing rules variations $(b)$ and $(c)$ we see that
the choice to have  face cards take the 
value $5$ rather than $10$ has a much larger effect on the failure
probability than the magician's choice of first key card position.
The final column of Table~5.2 presents the failure probabilities
for the i.i.d. uniform deck distribution on $\{ 1,2,...,10 \}$. One 
expects this  i.i.d.  uniform distribution to be comparable with rules 
variation $(c)$ rather than (b),
because the expected key card size is similar to case (c). (The Kruskal Count
$(c)$ mean value is slightly lower.)

To conclude: The rules variation to count face cards as having value $5$
rather than $10$ is important to the success of the Kruskal Count
trick in practice; the choice of the first card as key card offers a
further small improvement in success probability.

%
%

\section{ Appendix:  Proofs of Theorem 4.1, Lemma 4.2 and Theorem 4.2} 

{\bf Proof of Theorem 4.1.}
In view of Lemma 4.1, one has
\beql{eqA01}
M_N  \leq  {\rm Prob}[t < N]  \leq  NM_N 
\eeq
where
\beql{eqA02}
M_N  := \max_{1  \leq  k \leq  N} \left\{  \left(  1 - \df{1}{B} \right)^k  
~\mbox{Prob~$[t_{N,B}^{\ast \ast}   = k]$} \right\} ~.  
\eeq
The maximum will occur with $k  \approx  \gamma N$ for some constant
$\gamma  = \gamma (N,~B)$.
To estimate $\gamma$, we note that the travel time $t_{N,B}^{\ast \ast}$
beyond position N depends
on the successive transitions of the chain $\bar{\sL}_B$.
We convert this to a problem about successive states of the
{\em jump chain}
$\bar{\sL}_B^J$ having $2B(B  - 1)$ states which correspond to all
possible transitions of the chain $\bar{\sL}_B$.
A jump chain state $(i,~v)$ will mean state $i$ of $\bar{\sL}_B$
together with an allowable key card value $v$ which determines the next
state of $\bar{\sL}_B$.
The allowable values are
$1 \leq  v  \leq  B$ with $v \neq  |i|$.
The transition probability $P_{s,s'}$ from
$s  = (i,~v)$ to $s'  = (j,~v' )$ is $\frac{1}{B  - 1}$
when $j$ is uniquely determined by (\ref{eq401}) and
$1 \leq  v'  \leq  B$ with $v'  \neq  j$, and is 0 otherwise.

We let $  \{ ( i_k ,~ v_k  ) : k  + 1, 2, \cdot \cdot \cdot  \} $ denote a sequence
of states of $\bar{\sL}_B^J$, and introduce the modified travel time
\beql{eqA03}
\tilde{t}_{N,B} : = \min \{ k : v_1 + \ldots + v_k \geq 2N \} ~. 
\eeq
Now ${\rm Prob}[ \tilde{t}_{N,B}  \leq  \gamma N]$ can be estimated using
large deviation theory, using the following special case of
Theorem~1 of Donsker and Varadhan \cite{DonVar75}. \\
{\bf Theorem A.1.}  {\em For fixed $B$ and $\gamma$ one has 
as $N \rightarrow \infty$,}
\beql{eqA05}
{\rm Prob}[ \tilde{t}_{N,B}  \leq  \gamma N ]  = \exp ( - f_B ( \gamma ) N  +  o(N)) 
\eeq
{\em where the function $f_B ( \gamma )  = \gamma I^\ast ( \gamma )$ where}
\beql{eqA06}
I^\ast  ( \gamma )  := \inf \{  I( \mu ) : \mbox{weight} ( \mu ) \geq 
\df{2}{\gamma}  \}  ~. 
\eeq
{\em Here $\mu$ runs over the set of probability measures on the
state space $S$ of the chain $\bar{\sL}_B^J$, and}
$$
\mbox{ weight}( \mu )   := \sum_{s = (i,~v)  \in S} v  \mu ((i,~v)) ~,
$$
{\em is the expected card value for the measure $\mu$, and}
\beql{eqA07}
I( \mu )  :=  - \inf   \left\{  \sum_{s \in S}  \log 
\left( \df {\pi_u (s)}{u(s)} \right)  \mu (s) : u : 
S \rightarrow \RR^+ \right\}  ~, 
\eeq
{\em where $\pi_u (s)  := \sum_{s'  \in S} p_{s, s'} u (s' )$.}

It is easy to show that
$$
t_{N,B}^{\ast \ast}  \leq  \tilde{t}_{N,B} \leq  t_{N,B}^{\ast \ast} +  B~,
$$
and this yields
\beql{eqA04}
{\rm Prob}[ \tilde{t}_{N,B}  \leq  k ]  \leq  {\rm Prob}[t_{N,B}^{\ast \ast}  \leq  k]  \leq  {\rm Prob}
[ \tilde{t} _{N,B}  \leq  k  - B] ~.
\eeq
Using Theorem A.1 and (\ref{eqA04}) we see that the quantities
$t_{N,B}^{\ast \ast}$ and $\tilde{t}_{N,B}$ have the same
asymptotic behavior, with
\beql{eqA08}
{\rm Prob}[t_{N,B}^{\ast \ast}  \leq  \gamma N ]  =  \exp(  - f_B ( \gamma )
N  +  o (N)) 
\eeq
as $N  \rightarrow \infty$.
Now (\ref{eqA02}) leads us to define
\beql{eqA09}
 - \alpha_B  := \max_{0 < \gamma < \frac{B+1}{2}} 
\left[ \gamma  \log  \left(  1 -\df{1}{B} \right) - f_B ( \gamma ) \right]  ~.
\eeq
Using the fact that $f_B ( \gamma )$ is a strictly convex function, 
it can be shown
there is a unique value of  $\gamma$ attaining the maximum on the
right-hand side. With some further work this fact and
 (\ref{eqA01}), (\ref{eqA02}) and (\ref{eqA08}) imply that
$$
M_N  = \exp  ( - \alpha_B N  +  o(N))
$$
as $N \rightarrow \infty$.
The bound  of Theorem 4.1 follows.
$~~~\bsq$ \\

{\bf Proof of Lemma 4.2.}
We exhibit a $1-1$ correspondence between a sequence of states and 
admissible transitions for the
three chains, corresponding to moving two colored pebbles on the line
$\{n: n \geq 0 \}$ starting with both at zero.

We compare $\sL_B^-$ and $\bar{\sL}_B$.
The position after two pebble moves of $\sL_B^-$ corresponds to an
initial position for $\bar{\sL}_B$ together with the top key card and its
color.
Let $q_{ij}^-$, $q_{ij}$ denote the associated probability
distributions of the locations $(i,~j)$ of the white and black pebbles
for the chains $\bar{\sL}_B$, $\sL_B^-$ respectively,
which are: $q_{ij}^-  + \frac{1}{(B  - 1)^2}$ for
$1 \leq  i,~j \leq  B  - 1$ and 0 otherwise, $q_{ij}  +\frac{1}{B( B  - 1)}$
for $1 \leq  i,~j  \leq  B$ with $i \neq j$ and 0 otherwise. \\
{\bf Claim.} {\em There is a mapping $\phi^-$ of the probability mass $q_{ij}$
on $(i,~j)$ for $\bar{\sL}_B$ to various $(i' , j' )$
having $i'  \leq  i$ and $j'  \leq  j$ whose image is the
distribution $q_{ij}^-$.}

Let $(i,~j)  \succ  (i' ,~j' )$ mean
$\min (i,j) \geq \min (i', j')$ and
$\max (i,j) \geq \max (i', j')$, i.e. the
pebbles $(i,~j)$ are both moved further along the line then
$(i' ,~ j' )$, ignoring their colors.
Assuming the claim to be true for the moment, we have a
(stochastic) pairing of pebble positions
 $(i,~j)  \stackrel{\phi}{\rightarrow} (i' ,~ j' )$
such that $(i,~j)  \succ  (i' ,~ j' )$ between
$\bar{\sL}_B$ and $\sL_B^-$.
For each subsequent move, both claims have $B  - 1$ possible
transitions with probabilities each $\frac{1}{B  - 1}$.
For $\bar{\sL}_B$ in state $k  + i  - j$ the admissible value of the
next move is $ \{ 1,~2 , \ldots , B \}  -  \{ |k| \} $.
We map these transitions to transitions of $\sL_B^-$ in linear order,
with a mapping $\psi_{|k|}$ having
$\psi_{|k|} (i)  = i$ for $i  \leq  |k|$ and $\psi_j (i)  + i  - 1$
for $i \geq |k| + 1$.
One easily sees that if pebbles in $\bar{\sL}_B$ are at $(i_1 ,~ j_1 )$ and the corresponding
ones are at $(i_2 ,~ j_2 )$ with
$(i_1 ,~ j_1 )  \succ  (i_2 ,~ j_2 )$, and if the pebble
closer to the origin is moved $i$ resp. $\phi_{|k|} (i)$ for the two
claims resulting in positions $(i_1^\ast ,~ j_1^\ast )$,
$(i_2^\ast ,~ j_2^\ast )$ then $(i_1^\ast ,~ j_1^\ast )  \succ  (i_2^\ast ,~ j_2^\ast )$.
This gives a stochastic pairing of pebble positions at all subsequent
moves, with both pebbles of $\sL_B^-$ always being behind those of
$\bar{\sL}_B$ in the ordering $\succ$.
Consequently
\beql{eq412}
{\rm Prob} [t_{N,B}^{\ast \ast} > t ] \geq {\rm Prob}
[t_{N,B}^- > t ]
\eeq
for all $t$, which is the right side of (\ref{eq411}).

It remains to prove the claim.
Here we remark if $q_{ij} (n)$ and $q_{ij}^- (n)$ denote the
probability distribution of pebble locations of $\bar{\sL}_B$ and $\sL_B^-$ after $n$ moves then
$$
{\rm Prob} [t_{N,B}^{\ast \ast} > t ] = 1 -
\sum_{\stackrel{i \leq N}{j \leq N}} q_{ij} (t)
$$
with a similar formula for ${\rm Prob}[t_{N,B}^- >  t]$.
The proof above actually establishes the majorization inequalities
\beql{eq413}
\sum_{ \stackrel{i \leq i_0}{j \leq j_0}} q_{ij}^- (t) \geq
\sum_{ \stackrel{i \leq i_0}{j \leq j_0}} q_{ij} (t) ~,~~
\mbox{all ~ $i_0 \geq 1,~ j_0 \geq 1 $} ~,
\eeq
for all $t \geq 2$,
and the special case $i_0 = j_0 + N$
yields (\ref{eq412}).
The claim to be established is equivalent to proving that 
(\ref{eq413}) holds for the
case $t  = 2$.
Since the probabilities $q_{ij} = q_{ij} (2)$ and
$q'_{ij} = q'_{ij} (2)$ are explicitly known,
verifying (\ref{eq413}) is an easy
calculation.
The equivalence of the inequalities (\ref{eq413}) to the
 existence of a coordinate-monotone
probability rearrangement $\phi$ is a two-dimensional majorization inequality,
see Marshall and Olkin \cite{MarOlk79}.
(One can also prove the claim by explicitly constructing a suitable
mapping $\phi$ rather easily.)

The inequality (\ref{eq411}) relating $\bar{\sL}_B$ and $\sL_B^+$ 
is proved in
similar fashion.
If $q_{ij}^+  (t)$ is the probability that the pebbles
are at $(i,~j)$ after $t$ steps, then
\beql{eq414}
\sum_{\stackrel{i \leq i_0}{j \leq j_0}} q_{ij} (t) \geq
\sum_{\stackrel{i \leq i_0}{j \leq j_0}}
q_{ij}^+ (t) ~,~~ \mbox{all ~$ i_0 \geq 1 ~,~~ j_0 \geq 1 $}~,
\eeq
for all $t \geq 2$.
$~~~\bsq$ \\

{\bf Proof of Theorem 4.2.}
We let $A(N)  \approx  B(N)$ mean $A(N)  = B(N)^{1  +  o(1)}$ as
$N \rightarrow \infty$.

Consider first $P_{N,B}^-$.
Let $\tilde{t}_B^-(M)$ denote the travel time for the chain
$\sL_B^-$,  which counts the number of transitions up to and including
the transition
at which the sum of the jumps of the
chain exceeds $M$.
Then for any fixed sequence of transitions
\beql{eq419}
t_B^- (2N) \geq t_{N,B}^- \geq t_B^- (2N-N ) ~.
\eeq
Hence
\beql{eq420}
P_N^- \leq N \max_{1 \leq j \leq N} \left\{
\left( 1 - \df{1}{B} \right)^{j-1} ~\mbox{Prob$[ \tilde{t}_B^- (2N) + j ]$} \right\}
\eeq
and
\beql{eq420a}
P_N^- \geq \max_{1 \leq j \leq N} \left\{
\left( 1 - \df{1}{B} \right)^{j-1} ~\mbox{Prob$ [ \tilde{t}_B^- (2N-B) = j ]$} \right\} ~.
\eeq
It's easy to check that
\begin{eqnarray}
\label{eq421}
Q_N^- : =  \max_{1 \leq  j \leq  N}    \left\{   \left(  1 - \df{1}{B} \right)^{j-1}  {\rm Prob} [ \tilde{t}_B^- (2N)   = j]  \right\}   \nonumber \\
   \approx \max_{1 \leq j \leq N}
 \left\{  \left(  1 - \df{1}{B} \right)^{j-1} 
{\rm Prob} [ \tilde{t}_B^-
(2N - B)  + j ] \right\} 
\end{eqnarray}
using ${\rm Prob} [t_B^- (2N - B ) \geq j ] \geq {\rm Prob} [t_B^- (2N) \geq j + B ]$.
It suffices to asymptotically estimate $Q_N^-$.
One has
\beql{eq422}
Q_N^-  \approx   \sup  \left\{  B^{- \gamma N} \left( 
\df{\gamma N} {\gamma_1 N \gamma_2 N \cdot \cdot \cdot \gamma_{B-1} N} \right)  : 
\gamma_i \geq 0 ~,~~
\sum_{i=1}^{B-1} \gamma_i  = \gamma ~,~~
\sum_{i=1}^{B-1}~i \gamma_i  = 2 \right\}  ~.
\eeq
(Here $j  \approx \gamma N$ and $B^{- \gamma N}$ arises as
$\left( 1 - \frac{1}{B} \right)^{\gamma N} (B  - 1)^{- \gamma N}$.)
Using Stirling's formula, one obtains
$Q_N^-  \approx  \exp ( \frac {1}{2} \alpha_B^- N)$ where
$\frac {1}{2} \alpha_B^-$ is the optimal value of the constrained
maximization problem $(M^- )$ given by:
$$
\mbox{maximize} ~ Z  = - \gamma  \log  B  +  \gamma  \log  \gamma  - 
\sum_{i=1}^{B-1} \gamma_i  \log  \gamma_i
$$
\hspace*{.25in}{\em subject to}
\beql{eq423a}
\sum_{i=1}^{B-1} i \gamma_i   =   2~, 
\eeq
\beql{eq423b}
\sum_{i=1}^{B-1} \gamma_i  = \gamma ~, 
\eeq
\beql{eq423c}
\gamma_i \geq 0~ ~\mbox{for~ $1 \leq i \leq B - 1$} ~. 
\eeq
Introducing Lagrange multipliers $\lambda_1  ,~ \lambda_2$ for 
the two equality
constraints, and setting
$$
G  = \gamma  \log  B  +  \gamma  \log  \gamma - 
\sum_{i=1}^{B-1} \gamma_i  \log  \gamma_i - \lambda_1
\left(  \sum_{i=1}^{B-1} i \gamma_i - 2 \right)  - 
\lambda_2\left(  \sum_{i=1}^{B-1} \gamma_i - \gamma \right)~.
$$
Necessary conditions for an interior extremal are:
\beql{eq424a}
\df{\partial G}{\partial \gamma }   =   1  +  \log  \gamma - \log B  + 
\lambda_2 = 0~, 
\eeq
\beql{eq424b}
\df{\partial G}{\partial \gamma_i}   =  - 1 - \log  \gamma_i - i \lambda_1 - \lambda_2   = 0~. 
\eeq
These yield
\beql{eq425a}
\gamma   =  B \exp ( \lambda_2  - 1) ~, 
\eeq
\beql{eq425b}
\gamma_i  =  \exp ( \lambda_2  - 1)  \exp ( i \lambda_1 ) ~,~~~ 1 \le i \le B - 1. 
\eeq
Substituting these expressions into (\ref{eq423b})
and cancelling $\exp ( \lambda_2  - 1)$ from both sides yields
\setcounter{equation}{25}
\beql{eq426}
\sum_{i=1}^{B-1} \exp ( i \lambda_1 )   =  B~.
\eeq
Substituting the values above into  (\ref{eq424a}) and using this formula   yields
\beql{eq427}
\exp ( \lambda_2  - 1)  = \df{2}{\sum_{i=1}^{B-1}i  \exp ( i \lambda_1 )} 
\eeq
Now  (\ref{eq425a}) and (\ref{eq425b}) become
$$
\gamma   =   \df{2B}{\sum_{i=1}^{B-1} i \exp  ( i \lambda_1 )}
$$
$$
\gamma_i  = \df {2i  \exp ( i \lambda_1 )}
{\sum_{i=1}^{B-1} i \exp  ( i \lambda_1 )} ~.
$$
Using these formulas the objective function value $Z$ is evaluated
(with $F  = \sum_{i=1}^{B-1} i \exp  ( i \lambda_1 )$) as
\begin{eqnarray}
\label{eq428}
Z  & = &  \df{2B}{F} \left(  \log  \df{2B}{F} - \log  B \right)  - 
\sum_{i=1}^{B-1}\df{2  \exp (i \lambda_1 )}{F}  \log 
\left( \df {2  \exp  ( i \lambda_1 )}{F} \right) \nonumber \\
 & = & \df{2B}{F}  \log  \df{2}{F} -  \df{2}{F} \sum_{i=1}^{B-1} \exp ( i \lambda_1 ) \left[  \log \df{2}{F}  +  i \lambda_1 \right]  \nonumber \\
& = & \df{2 \lambda_1}{F}  \sum_{i=1}^{B-1} i  \exp  (i \lambda_1 )  =  - 2 \lambda_1 ~.
\end{eqnarray}
Combining this with (\ref{eq426}) gives (\ref{eq416}) for $P_{N,B}^-$ 
with $\frac {1}{2} \alpha_B^-$ given by ({\ref{eq417a}),
provided the maximum of $(M^- )$ occurs
at an interior point where all $\gamma_i >  0$.
We omit the details of checking that boundary
extremals having some $\gamma_i  = 0$ do not give the absolute
maximum in $(M^- )$.

The case of $P_{N,B}^+$ is handled by analogous arguments.
One reduces it to solving the 
constrained maximization problem $(M^+)$ given by:
$$
 \mbox{maximize} ~~ Z =  - \gamma  \log B  +  \gamma  \log 
\gamma - \sum_{i=1}^{B-1} \gamma_i  \log  \gamma_i
$$
\hspace*{.5in}{\em subject to}
\beql{eq429a}
\sum_{i=1}^{B-1} (i  +  1) \gamma_i  = 2~, 
\eeq
\beql{eq429b}
\sum_{i=1}^{B-1} \gamma_i  = \gamma ~,  
\eeq
\beql{eq429c}
\gamma_i \geq 0 ~~\mbox{for ~$ 1 \leq i \leq B$} ~. 
\eeq
Again $Z  = - 2  \lambda_1$ at the extremal point, and $\lambda_1 $
is determined by
$$
\sum_{i=1}^{B-1} \exp(i  +  1) \lambda_i  = B~,
$$
which is (\ref{eq417b}).

The asymptotic formulae ({\ref{eq418a}) and (\ref{eq418b}) are obtained
from the formulas (\ref{eq417a}) and (\ref{eq417b}),
by setting 
$\exp (\frac {1}{2} \alpha )   = 1  +  \frac{2}{B^2}  + 
\frac{\delta}{B^3}  +  O (B^{-4})$ , in which $\delta$ is an
unknown to be determined. We find it by noting that
$$
\exp (\frac {1}{2} i \alpha )  = 1  +  i \left(  \df{2}{B^2 }  + \df{\delta}{B^3 } \right)  + 
i^2   \df{2}{B^4}  +  O \left(  B^{-3} \left(  1  +  \df{i}{B} \right)^{-1}\right ) ~,
$$
 and solving for the appropriate value of $\delta$.~~ \bsq

{\tt
\begin{tabular}{ll}
email: & jcl@research.att.com \\
 &  rains@research.att.com \\
& rvdb@princeton.edu
\end{tabular}
 }
\end{document}